\documentstyle[11pt,A4]{article}
\def\QED{{\hfill$\Box$}}
\begin{document}
\bibliographystyle{plain}
\title{On computing factors of cyclotomic polynomials%
\thanks{Copyright \copyright\ 1992--2010, R.~P.~Brent.
	\hfill rpb135 typeset using \LaTeX}}
\author{Richard P.\ Brent\\
	Computer Sciences Laboratory\\
	Australian National University\\
	Canberra, ACT 0200\\
	\\
	Report TR-CS-92-13\\
	September 1992\\
	\\
	{\em In memory of}\\
	{\em Derrick H.~Lehmer}\\
	{\em 1905--1991}}
\date{~}
\maketitle
\thispagestyle{empty}			%

\vspace{-1cm}
\begin{abstract}
For odd square-free $n > 1$ the cyclotomic polynomial $\Phi_n(x)$
satisfies the identity of Gauss %
\[	4\Phi_n(x) = A_n^2 - (-1)^{(n-1)/2}nB_n^2 .	\]
A similar identity of Aurifeuille, Le Lasseur and Lucas is
\[	\Phi_n((-1)^{(n-1)/2}x) = C_n^2 - nxD_n^2	\]
or, in the case that $n$ is even and square-free,
\[	\pm\Phi_{n/2}(-x^2) = C_n^2 - nxD_n^2,	\]
Here $A_n(x),\ldots, D_n(x)$ are polynomials with integer coefficients.
We show how these coefficients can be computed
by simple algorithms which require $O(n^2)$ arithmetic operations
and work over the integers.
We also give explicit formulae and generating functions
for $A_n(x),\ldots, D_n(x)$, and illustrate the application to integer
factorization with some numerical examples.
\smallskip

\smallskip


\end{abstract}

\vfill\eject		%
\section{Introduction}

For integer $n > 0$ let $\Phi_n(x)$ denote the
cyclotomic polynomial
\begin{equation}
  \Phi_n(x) = \prod_{\scriptstyle 0 < j \le n \atop	%
	\scriptstyle  (j,n)=1} (x - \zeta^j), 		\label{eq:C1}
\end{equation}
where $\zeta$ is a primitive $n$-th root of unity.
Clearly
\begin{equation}
  x^n - 1 = \prod_{d \vert n} \Phi_d(x) ,		\label{eq:C2}
\end{equation}
and the M\"obius inversion formula~\cite{Har60} gives
\begin{equation}
  \Phi_n(x) = \prod_{d \vert n} (x^d - 1)^{\mu(n/d)} .	\label{eq:C3}
\end{equation}

Equation~(\ref{eq:C1}) is useful for theoretical purposes,
but~(\ref{eq:C3}) is more convenient for computation as it leads to a
simple algorithm for computing the coefficients of $\Phi_n(x)$
or evaluating $\Phi_n(x)$ at integer arguments using only integer arithmetic.
If $n$ is square-free the relations
\begin{equation}
\Phi_n(x) = \cases{	%
	x-1,				&if $n = 1$;\cr
	\Phi_{n/p}(x^p)/\Phi_{n/p}(x),	&if $p \vert n$,
					 $p$ prime;\cr}
\label{eq:phirec}
\end{equation}
give another convenient recursion for computing $\Phi_n(x)$.

Although $\Phi_n(x)$ is irreducible over $Z$ (see for example~\cite{Wae53}),
$\Phi_n(x)$ may be reducible over certain quadratic fields.
For example,
\begin{equation}
	4\Phi_5(x) = (2x^2 + x + 2)^2 - 5x^2,
\end{equation}
so $\Phi_5(x)$ has factors $x^2 + \left({{1 \pm \sqrt{5}}\over 2}\right)x + 1$
whose coefficients are algebraic integers in $Q[\sqrt{5}]$.

For odd square-free $n > 1$ the cyclotomic polynomial $\Phi_n(x)$
satisfies the identity
\begin{equation}
  4\Phi_n(x) = A_n^2 - (-1)^{(n-1)/2}nB_n^2 .		\label{eq:G1}
\end{equation}
Gauss~\cite{Gau01} proved~(\ref{eq:G1}) for odd prime $n$;
the generalization to other odd square-free $n$ is due to
Dirichlet~\cite{Dir94}.
Related identities of Aurifeuille and Le Lasseur~\cite{Aur77} are
\begin{equation}
	\Phi_n((-1)^{(n-1)/2}x) = C_n^2 - nxD_n^2       \label{eq:A1}
\end{equation}
for odd square-free $n$, and
\begin{equation}
	\Phi_{n/2}(-x^2) = C_n^2 - nxD_n^2		\label{eq:A2}
\end{equation}
for even square-free $n > 2$.	%
For a proof, see Lucas~\cite{Luc78b} or Schinzel~\cite{Sch62}.

In~(\ref{eq:G1}--\ref{eq:A2}),
$A_n(x), \ldots , D_n(x)$ are polynomials with integer coefficients,
and without loss of generality we can assume that
$A_n(x)/2$, $B_n(x)$, $C_n(x)$ and $D_n(x)$ are monic.
In Section~\ref{sec:Algs}
we show how the coefficients of $A_n, \ldots, D_n$ can be computed
by simple algorithms which require $O(n^2)$ arithmetic operations
and work entirely over the integers.

In Section~\ref{subsec:no} we summarize our notation for future reference.
Some numerical examples are given in
Sections~\ref{subsec:ex}--\ref{subsec:Beeger}, and
Newton's identities are discussed in Section~\ref{subsec:Newton}.
Then, in Section~\ref{sec:Theory}, we discuss the theoretical basis for
the algorithms. The results for $A_n$ and $B_n$ are known (though
perhaps forgotten)~-- they may be found in Dirichlet~\cite{Dir94}.
We present them in Section~\ref{subsec:Gauss}
for the sake of completeness and to aid the reader in
understanding the results for $C_n$ and $D_n$.

The algorithms are presented in Section~\ref{sec:Algs}.
The algorithm (Algorithm~D)
for computing $A_n$ and $B_n$ is essentially due to
Dirichlet~\cite{Dir94}, who illustrated it with
some numerical examples but did not state it in general terms.
The algorithm (Algorithm~L) for computing $C_n$ and $D_n$ appears
to be new.
In Section~\ref{subsec:Stev}
we comment briefly on Stevenhagen's algorithm~\cite{Ste87}
and compare it with Algorithm~L.

Finally, in Section~\ref{sec:GF} we give some explicit formulas
for $A_n(x), \ldots, D_n(x)$. These may be regarded as generating
functions if $x$ is an indeterminate, or may be used to compute
$A_n(x), \ldots, D_n(x)$ for given argument $x$.
In the special case $x = 1$ the results for $A_n(1), B_n(1)$ reduce to
known formulas involving the
class number of the quadratic field
$Q[\sqrt{\pm n}]$.                	%

One application of cyclotomic polynomials
is to the factorization of
integers of the form $a^n \pm b^n$: see for
example~\cite{Bre92,Bri83,Cun15,Cun25,Kra22,Kra24,Rie85,Sch62,Ste87}.
If $x = m^2n$ for any integer $m$, then~(\ref{eq:A1}--\ref{eq:A2})
are differences of squares, giving
rational integer factors of $x^n \pm 1$.
Examples may be found in Section~\ref{subsec:intfac}.
For the reader interested in integer factorization, our most
significant results are Algorithm~L of Section~\ref{subsec:AlgsCD}
and Theorem~\ref{thm:Ta} of Section~\ref{subsec:intfac}.

\subsection{Notation}
\label{subsec:no}

Unless qualified by ``algebraic'', the term ``integer'' means
a rational integer.
$x$ usually denotes an indeterminate,
occasionally a real or complex variable.

$\mu(n)$ denotes the M\"obius function,
$\phi(n)$ denotes Euler's totient function, and
$(m,n)$ denotes the greatest common divisor of $m$ and $n$.
For definitions and properties of these functions,
see for example~\cite{Har60}. Note that $\mu(1) = \phi(1) = 1$.

$(m \vert n)$ %
denotes the Jacobi symbol\footnote{See, for example, Riesel~\cite{Rie85}.
To avoid ambiguity, we {\em never} write the Jacobi symbol as
$\left( m \over n \right)$.
Note that $m \vert n$ without parentheses means that $m$ divides $n$.}
except that, as is usual for the Kronecker symbol\footnote{See, for example,
Landau~\cite{Lan27}.},
$(m \vert n)$ is defined as 0 if $(m,n) > 1$.
Thus, when specifying a condition such as $(m \vert n) = 1$
we may omit the condition $(m,n) = 1$.

$n$ denotes a positive integer (square-free from
Section~\ref{subsec:Gauss} on).
For given $n$, we define integers $n'$, $s$ and $s'$ as follows:
\[n' = \cases{n, 	&if $n = 1 \bmod 4$;\cr
	     2n, 	&otherwise.\cr} \]
\[s  = \cases{-1,	&if $n = 3 \bmod 4$;\cr
	     +1,	&otherwise.\cr}	\]
\[s' = \cases{-1,	&if $n = 5 \bmod 8$;\cr
	      +1,	&otherwise.\cr} \]
It is convenient to write $g_k$ for $(k,n)$ and $g_k'$ for $(k,n')$.

Define
\begin{equation}
  F_n(x) = \cases{\Phi_n(sx),				&if $n$ is odd;\cr
		  (-1)^{\phi(n/2)}\Phi_{n/2}(-x^2),	&if $n$ is even.\cr}
							\label{eq:F_n_def}
\end{equation}
Thus
we can write~(\ref{eq:A1}--\ref{eq:A2}) as
\begin{equation}
	F_n(x) = C_n^2 - nxD_n^2 .			\label{eq:Fn2}
\end{equation}
The factor $(-1)^{\phi(n/2)}$ in the definition of $F_n$ is only relevant
if $n = 2$, and ensures that~(\ref{eq:Fn2}) is valid
for $n = 2$ (with $C_2(x) = x+1$, $D_2(x) = 1$).
The Aurifeuillian factors of $F_n(x)$ are
\[ F_n^{+}(x) = C_n(x) + \sqrt{nx}D_n(x)	 \]
and
\[ F_n^{-}(x) = C_n(x) - \sqrt{nx}D_n(x). \]
From~(\ref{eq:Fn2}) we have $F_n(x) = F_n^{-}(x)F_n^{+}(x)$.
We may write $F_n^{\pm}$ for one of $F_n^{+}$, $F_n^{-}$.

We sometimes need to specify a particular complex square root.
If $m < 0$ then $\sqrt{m}$ means $i\sqrt{\vert m \vert}$.

$d$ is usually the degree of a polynomial,
while $D$ is the discriminant of a quadratic form.
For odd square-free $n$ we always have $D = sn$,
so $D = 1 \bmod 4$.

Some additional notation is introduced in Section~\ref{subsec:Newton}.

\subsection{Examples}
\label{subsec:ex}

Taking $n = 15$, we have
\[ \Phi_{15}(x) = {(x^{15}-1)(x-1) \over (x^5-1)(x^3-1)}
	= x^8 - x^7 + x^5 - x^4 + x^3 - x + 1 ,	\]
\[ A_{15}(x) = 2x^4 - x^3 - 4x^2 - x + 2 ,	\]
\[ B_{15}(x) = x^3 - x , \]
\[ C_{15}(x) = x^4 + 8x^3 + 13x^2 + 8x + 1 ,	\]
\[ D_{15}(x) = x^3 + 3x^2 + 3x + 1 ,		\]
and the reader may easily verify that~(\ref{eq:G1}) and~(\ref{eq:A1})
are satisfied.  As an example of~(\ref{eq:A2}), for $n = 14$ we have
\[ F_{14}(x) = \Phi_7(-x^2) = {x^{14}+1 \over x^2+1}
	= x^{12} - x^{10} + x^8 - x^6 + x^4 - x^2 + 1, \]
\[ C_{14}(x) = x^6 + 7x^5 + 3x^4 - 7x^3 + 3x^2 + 7x + 1, \]
and
\[ D_{14}(x) = x^5 + 2x^4 - x^3 - x^2 + 2x + 1 .	\]

\subsection{The identities of Beeger and Schinzel}
\label{subsec:Beeger}

Taking $n = 5$ in~(\ref{eq:A1}) we obtain
\[ \Phi_5(x) = (x^2 + 3x + 1)^2 - 5x(x+1)^2 , \]
so $\Phi_5(x^2)$ has factors $x^4 + 3x^2 + 1 \pm \sqrt{5}(x^3 + x)$
in $Q[\sqrt{5}]$.  Replacing $x$ by $x^3$ we obtain factors of
$\Phi_5(x^6)$.  Now
\[ \Phi_{15}(x^2) = \Phi_5(x^6)/\Phi_5(x^2) ,	\]
and by division (taking the factors with opposite signs of $\sqrt{5}$)
we obtain factors
\[	x^8 + 2x^6 + 3x^4 + 2x^2 + 1 \pm \sqrt{5}(x^7 + x^5 + x^3 + x) \]
of $\Phi_{15}(x^2)$.  Thus
\begin{equation}
 \Phi_{15}(x) = (x^4 + 2x^3 + 3x^2 + 2x + 1)^2
			- 5x(x^3 + x^2 + x + 1)^2 .	\label{eq:Sch1}
\end{equation}
This is not of the form~(\ref{eq:Fn2}) because it gives a
factorization of $\Phi_{15}(\pm x)$ over
$Q[\sqrt{\pm 5}]$ instead of $Q[\sqrt{\mp 15}]$.	%
Instead,~(\ref{eq:Sch1}) is an example of the more general identities
of Beeger~\cite{Bee51} and Schinzel~\cite{Sch62}.
These identities can all be obtained in a similar manner
from~(\ref{eq:Fn2}), so in the application to
integer factorization they do not give any factors which could not
be found from several applications of~(\ref{eq:Fn2})
and some greatest common divisor calculations.
For this reason we have restricted our attention
to identities of the form~(\ref{eq:G1}) and~(\ref{eq:Fn2}).

\subsection{Newton's identities}
\label{subsec:Newton}

Let
\[	P(x) = \prod_{j=1}^d (x - \xi_j) = \sum_{j=0}^d a_jx^{d-j}	\]
be a polynomial of degree $d$ with arbitrary roots $\xi_j$ and coefficients
$a_0 = 1, a_1, \ldots, a_d$.

For $k > 0$, define
\begin{equation}
	p_k = \sum_{j=1}^d \xi_j^k .		\label{eq:pk}
\end{equation}
Newton (1707)\footnote{See Turnbull~\cite{Tur52},
where the notation $s_k$ is used in place of our $p_k$.
It would be confusing to use Turnbull's notation because
we have used $s$ and $s'$ for other purposes.
Note that our $p_k$ are {\em not} generally prime numbers. %
}
showed how to express the elementary symmetric functions
$a_1, a_2, \ldots$ in terms of the sums of powers $p_1, p_2, \ldots$

We may find $a_1, \ldots, a_d$
by solving a lower triangular linear system of
special form~\cite{Tur52}. %
Writing the solution explicitly in the form of a linear
recurrence, we have
\begin{equation}
 ka_k = - \sum_{j=0}^{k-1} p_{k-j}a_j  \label{eq:NID}
\end{equation}
for $k = 1, \ldots, d$.
An alternative expression for $a_k$ as a %
determinant may be obtained by applying Cramer's rule to the
lower triangular system.
However, for computational purposes~(\ref{eq:NID}) is more convenient.

In Section~\ref{sec:GF} we use the following generating function~\cite{Rio78}
for $(a_0,a_1,\ldots)$:
\begin{equation}
	x^dP(1/x) = \sum_{j=0}^d a_jx^j =
	\exp\left(-\sum_{j=1}^\infty p_jx^j/j\right)		\label{eq:GF1}
\end{equation}
Differentiating both sides of~(\ref{eq:GF1}) and equating coefficients shows
that~(\ref{eq:NID}) and~(\ref{eq:GF1}) are formally equivalent.
An independent proof of~(\ref{eq:GF1}) is the following:
for sufficiently small $x$ we have	%
\[ \ln(x^dP(1/x)) = \sum_{k=1}^d \ln(1 - \xi_k x)
	= -\sum_{k=1}^d \sum_{j=1}^{\infty} \xi_k^j x^j/j
	= -\sum_{j=1}^{\infty} \left(\sum_{k=1}^d \xi_k^j\right) x^j/j
	= -\sum_{j=1}^{\infty} p_jx^j/j . \]
In all our applications of~(\ref{eq:GF1})
the $p_j$ are bounded, so the infinite series
converges for $\vert x \vert < 1$.

In the following, $p_j$ and $a_j$ are not fixed, but depend on the
particular polynomial under consideration at the time.  This should
not cause any confusion.

\section{Theoretical Basis for the Algorithms}
\label{sec:Theory}

Our idea
is to compute sums of powers of
certain	%
roots of the polynomials
occurring on the left side of~(\ref{eq:G1}) and~(\ref{eq:Fn2}), and then
use Newton's identities in the form~(\ref{eq:NID}) to compute the
coefficients of $A_n, \ldots, D_n$.

\subsection{Cyclotomic polynomials}
\label{subsec:cyclo}

First consider the computation of
the coefficients of the cyclotomic polynomial $\Phi_n(x)$ for $n > 1$.
This is presented to illustrate a simple case of the technique;
in practice it is more efficient to compute $\Phi_n(x)$ from~(\ref{eq:C3}).

Let $\zeta$ be a primitive $n$-th root of unity.  To apply Newton's
identities we need to evaluate
\begin{equation}
 	p_k = \sum_{\scriptstyle 0 < j < n \atop
	\scriptstyle  (j,n)=1} \zeta^{jk} 		\label{eq:C4}
\end{equation}
for $k = 1, 2, \ldots, \phi(n)$.
This problem is well-known\footnote{If
$\zeta = e^{2{\pi}i/n}$ then our $p_k$ is ``Ramanujan's sum'' $c_n(k)$,
in the notation of
Ramanujan~\cite{Ram18}	%
or Chapter 26 of Davenport~\cite{Dav80}.}.

If $n$ is prime the problem is easy: from
\[ 1 + \zeta + \zeta^2 + \cdots + \zeta^{n-1} =
	{{1 - \zeta^n}\over{1 - \zeta}} = 0 ,	\]
we have $p_1 = -1$.  Moreover, for any $k$ with $(k,n) = 1$,
the map $z \mapsto z^k$ merely permutes
$\{\zeta, \ldots, \zeta^{n-1}\}$, so $p_k = p_1$.

Now consider the general case, $n$ not necessarily prime.
From~(\ref{eq:C3}) it is clear that $p_1 = \mu(n)$.
Let $g_k = (k,n)$.  If $g_k = 1$ then the same argument as before
shows that $p_k = p_1$.  If $g_k > 1$ then the sum~(\ref{eq:C4})
defining $p_k$ consists of $\phi(n)/\phi(n/g_k)$ copies of a sum of
primitive $(n/g_k)$-th roots of unity.
Thus, the result is		%
\begin{equation}
	p_k = {{\mu(n/g_k)\phi(n)} \over \phi(n/g_k)} .	\label{eq:C5}
\end{equation}
Using~(\ref{eq:C5}) the coefficients
$a_1, \ldots, a_{\phi(n)}$ of $\Phi_n(x)$ may be
evaluated from the recurrence~(\ref{eq:NID}).

As an application of~(\ref{eq:GF1}) and~(\ref{eq:C5}) we prove two Lemmas
which give upper bounds on $\vert \Phi_n(x) \vert$
and $\vert F_n(x) \vert$ for complex $x$ outside the unit circle.
Here $F_n(x)$ is the modified cyclotomic polynomial
defined by~(\ref{eq:F_n_def}).
Lemma~\ref{lemma:Lb} is used in Section~\ref{subsec:intfac}.

\newtheorem{lemmaa}{Lemma}
\begin{lemmaa}
If $\vert x \vert \ge R > 1$ then
\[ \vert \Phi_n(x) \vert < R^{\phi(n)}\exp\left({1 \over R-1}\right) .\]
\label{lemma:La}
\end{lemmaa}

\leftline{\bf Proof}
\medskip
Let $d = \phi(n)$.
From~(\ref{eq:GF1}) with $x$ replaced by $1/x$, we have
\[ \Phi_n(x)/x^d = \exp\left(-\sum_{j=1}^\infty p_jx^{-j}/j\right) ,\]
but from~(\ref{eq:C5}) we have
\[ \vert p_j \vert \le g_j \le \min(j,n) ,	\]
so
\[ \vert \Phi_n(x)/x^d \vert < \exp\left(\sum_{j=1}^\infty R^{-j} \right)
	= \exp\left({1 \over R-1}\right) .	\]
This completes the proof.					\QED
\medskip

\newtheorem{lemmab}[lemmaa]{Lemma}
\begin{lemmab}
If $n > 1$ is square-free and $\vert x \vert \ge R > 1$ then
\[ \vert F_n(x) \vert < R^{\phi(2n)}\exp\left({1 \over R-1}\right) . \]
\label{lemma:Lb}
\end{lemmab}

\leftline{\bf Proof}
\medskip
From the definition~(\ref{eq:F_n_def}) of $F_n(x)$,
\[ {\rm deg\;} F_n = \cases{\phi(n),	&if $n$ is odd;\cr
			    2\phi(n/2),	&if $n$ is even;\cr}	\]
so it is easy to see that
\[ {\rm deg\;} F_n = \phi(2n)	\]
in both cases.
The bound on $\vert F_n(x) \vert$ follows from Lemma~\ref{lemma:La}
applied to $\Phi_n(\pm x)$ if $n$ is odd,
and from Lemma~\ref{lemma:La}
applied to $\Phi_{n/2}(-x^2)$ if $n$ is even.			\QED

\subsection{The identity of Gauss}
\label{subsec:Gauss}

From now on we assume that $n > 1$ is square-free.
In this subsection we also assume that $n$ is odd.
Consider the polynomial
\begin{equation}
  G_n(x) = \prod_{\scriptstyle 0 < j < n \atop
	\scriptstyle  (j \vert n)=1} (x - \zeta^j) 		\label{eq:G2}
\end{equation}
of degree $\phi(n)/2$,
where $\zeta = e^{2{\pi}i/n}$ is a primitive $n$-th root of unity.
The particular choice of primitive root is only significant for the
sign of the square root of $sn$ appearing in the equations below.

From Dirichlet~\cite{Dir94},
\begin{equation}
	2G_n(x) = A_n(x) - \sqrt{sn}B_n(x) ,        		\label{eq:G2.5}
\end{equation}
where $A_n$ and $B_n$ are as in~(\ref{eq:G1}),
and $s$ is defined in Section~\ref{subsec:no}.
Since $n$ is odd, we have $s = (-1 \vert n) = (-1)^{(n-1)/2}$.

Define
\begin{equation}
	2{\tilde G}_n(x) = A_n(x) + \sqrt{sn}B_n(x) ,        \label{eq:Gtilde}
\end{equation}
so Gauss's identity~(\ref{eq:G1}) may be written as
\begin{equation}
	\Phi_n(x) = G_n(x){\tilde G}_n(x) .			\label{eq:G2.6}
\end{equation}

The sums of $k$-th powers of roots of $G_n(x)$ are
\begin{equation}
  	p_k = \sum_{\scriptstyle 0 < j < n \atop
	\scriptstyle  (j \vert n)=1} \zeta^{jk} . 		\label{eq:G3}
\end{equation}
Let $g_k = (k,n)$. Then
\begin{equation}
   2p_k = \cases{ \mu(n) + (k \vert n)\sqrt{sn}, &if $g_k = 1$;\cr
                  \mu(n/g_k)\phi(g_k), &otherwise.\cr}		\label{eq:G5}
\end{equation}
The result~(\ref{eq:G5}) is essentially due to Dirichlet~\cite{Dir94},
but we sketch a proof.
If $g_k = 1$, then~(\ref{eq:G5}) follows
from the discussion in Section~\ref{subsec:cyclo}
(where $p_k$ has a different meaning!)
and the classical result	%
that the Gaussian sum
\[	\sum_{0 < j < n} (j \vert n)\zeta^j \]	%
is $\sqrt{sn}$.
On the other hand, if $g_k > 1$, observe that $(g_k, n/g_k) = 1$
because $n$ is square-free.
Thus, we can write the summation index $j$ in~(\ref{eq:G3}) in the form
$j = j_0g_k + j_1(n/g_k)$,
and $(j \vert n) = (j_1 \vert g_k)(j_0 \vert (n/g_k))$.
Since $jk = j_0{g_k}k \bmod n$, $\zeta^{jk}$ is independent of $j_1$,
and it follows that the sum~(\ref{eq:G3}) defining $p_k$
consists of $\phi(g_k)/2$ copies of a complete sum of
primitive $(n/g_k)$-th roots of unity.
Thus,~(\ref{eq:G5}) follows as in the proof
of~(\ref{eq:C5}).	   %

Although (\ref{eq:G5}) has been written with two cases for the sake of
clarity, our convention that $(k \vert n) = 0$ if $(k,n) > 1$ implies that
the expression
\begin{equation}
   2p_k = \mu(n/g_k)\phi(g_k) + (k \vert n)\sqrt{sn}	\label{eq:G6}
\end{equation}
is valid in both cases.
Similarly for ${\tilde G}_n(x)$,
with the sign of $\sqrt{sn}$ in~(\ref{eq:G6}) reversed.

Observe that $p_k \in Q[\sqrt{n}]$ is real if $n = 1 \bmod 4$,
but $p_k \in Q[\sqrt{-n}]$ is complex if $g_k = 1$ and $n = 3 \bmod 4$.

Using~(\ref{eq:G5}), the coefficients of $G_n(x)$, and hence of
$A_n(x)$ and $B_n(x)$, may be evaluated from the recurrence~(\ref{eq:NID}).
Moreover, it is possible to perform the computation using only integer
arithmetic.
Details are given in Section~\ref{sec:Algs}.

\subsection{The identities of Aurifeuille, Le Lasseur and Lucas}
\label{subsec:Lucas}

Here we assume that $n > 1$ is square-free, but not necessarily odd.
Recall the definitions of $n'$, $s$ and $s'$ from Section~\ref{subsec:no}.

Let $\zeta = e^{{\pi}i/n'}$ be a primitive $2n'$-th root of unity.
The particular choice of primitive root is only significant
for the sign of the square root in~(\ref{eq:A5}).
Consider the polynomial
\begin{equation}
  L_n(x) = \prod_{j \in S_n} (x - \zeta^j) ,			\label{eq:A3}
\end{equation}
where
\begin{equation}
  S_n = \cases{
	\{j \;\vert\; 0 < j < 2n',\; (j,n') = 1,\; (j \vert n)=(-1)^j\},
		&if $n = 1 \bmod 4$;\cr
		&~\cr				%
        \{j \;\vert\; 0 < j < 2n',\; (j,n') = 1,\; (n \vert j)=1\},
		&otherwise.\cr}
								\label{eq:A4}
\end{equation}

Observe that $L_n(x)$ has degree $\phi(n') = \phi(2n)$.
Also, $j \in S_n$ iff $2n'-j \in S_n$, so the coefficients of $L_n(x)$ are real.
In fact, from~(\ref{eq:A5}) below, they are in $Q[\sqrt{n}]$.
We later use the fact that $L_n(x)$ is symmetric.

Schinzel~\cite{Sch62} essentially shows (with a different notation) that
\begin{equation}
	L_n(x) = C_n(x^2) - s'x\sqrt{n}D_n(x^2)			\label{eq:A4.5}
\end{equation}
where $C_n(x)$ and $D_n(x)$ are the polynomials of~(\ref{eq:Fn2}).
Define
\begin{equation}
	{\tilde L}_n(x) = L_n(-x) = C_n(x^2) + s'x\sqrt{n}D_n(x^2) ,    
							     \label{eq:Ltilde}
\end{equation}
so after a change of variable~(\ref{eq:Fn2}) may be written as
\begin{equation}
	F_n(x^2) = L_n(x){\tilde L}_n(x) .			\label{eq:A4.7}
\end{equation}
Clearly $F_n^{-}(x) = L_n(s'\sqrt{x})$
and $F_n^{+}(x) = {\tilde L}_n(s'\sqrt{x})$.

Let $g_k' = (k,n')$.
The sums $p_k$ of $k$-th powers of roots of $L_n(x)$ are

\begin{equation}
  p_k = \cases{	(n \vert k)s'\sqrt{n},			   &if $k$ is odd;\cr
	\mu(n'/g_k')\phi(g_k')\cos\left((n-1)k\pi/4\right),&if $k$ is even.\cr}
								\label{eq:A5}
\end{equation}
Observe that the cosine in~(\ref{eq:A5}) is 0 or $\pm 1$, and
depends only on $n \bmod 4$ and $k/2 \bmod 4$.
The proof of~(\ref{eq:A5}) is similar to that of~(\ref{eq:G5}),
but tedious because of the number of
cases to be considered. Thus, we omit the details.

Using~(\ref{eq:A5}), the coefficients of $L_n(x)$, and hence of
$C_n(x)$ and $D_n(x)$, may be evaluated from the recurrence~(\ref{eq:NID}).
Details are given in Section~\ref{sec:Algs}.

\section{Algorithms}
\label{sec:Algs}

In this section we use the analytic results of Section~\ref{sec:Theory}
to derive efficient
algorithms for computing the coefficients of the polynomials
$A_n, \ldots, D_n$.

\subsection{An algorithm for computing $A_n$ and $B_n$}
\label{subsec:AlgsAB}

Consider the computation of $A_n$ and $B_n$ for
odd square-free $n$. Our notation is the same as in
Section~\ref{subsec:Gauss}.
Write
\[	A_n(x) = \sum_{j=0}^d \alpha_jx^{d-j} ,	\]
\[	B_n(x) = \sum_{j=0}^d \beta_jx^{d-j} ,	\]
where $d = \phi(n)/2$, $\alpha_0 = 2$, $\beta_0 = 0$, $\beta_1 = 1$.

Recall the definition~(\ref{eq:G3}) of $p_k$.
For $k > 0$ we have, from~(\ref{eq:G6}),
\[	2p_k = q_k + r_k\sqrt{sn} ,	\]
where $q_k$ and $r_k$ are integers given by
\begin{equation}
	q_k = \mu(n/g_k)\phi(g_k)				\label{eq:AL1}
\end{equation}
and
\begin{equation}
	r_k = (k \vert n) .					\label{eq:AL2}
\end{equation}
Using~(\ref{eq:NID}) and~(\ref{eq:G2.5}), we obtain the recurrences
\begin{equation}
	\alpha_k = {1 \over 2k}\sum_{j=0}^{k-1}
		\left(snr_{k-j}\beta_j - q_{k-j}\alpha_j\right)	\label{eq:AL3}
\end{equation}
and
\begin{equation}
	\beta_k  = {1 \over 2k}\sum_{j=0}^{k-1}
		\left(r_{k-j}\alpha_j - q_{k-j}\beta_j\right)	\label{eq:AL4}
\end{equation}
for $k = 1, 2, \ldots, d$.
The algorithm is now clear:
\medskip

\leftline{{\bf Algorithm D} (for Dirichlet)}

\nopagebreak
\begin{itemize}
\item[1.] Evaluate $q_k$ and $r_k$ for $k = 1, \ldots, d$
	  using~(\ref{eq:AL1}--\ref{eq:AL2}).
\item[2.] Set $\alpha_0 \leftarrow 2$ and $\beta_0 \leftarrow 0$.
\item[3.] Evaluate $\alpha_k$ and $\beta_k$ for $k = 1, \ldots, d$
	  using~(\ref{eq:AL3}--\ref{eq:AL4}).
\end{itemize}
\medskip

\leftline{\bf Comments on Algorithm D}

\nopagebreak
\begin{itemize}
\item[1.] (\ref{eq:AL3}--\ref{eq:AL4}) should give exact integer results;
in practice a sum not divisible by $2k$ is a symptom of integer overflow.
\item[2.] The operation count can be reduced by a factor of close to four if
advantage is taken of the following properties of $A_n$ and $B_n$:
\item[]   $A_n$ is anti-symmetric if its degree $d = \phi(n)/2$ is odd,
	  otherwise $A_n$ is symmetric (except for the trivial case
	  $A_3(x) = 2x+1$).
	  Thus, we may use
\[	\alpha_k = (-1)^d\alpha_{d-k} 	\]
	if $2k > d$ and $n > 3$.
\item[] $B_n/x$ is antisymmetric if $n$ is composite and $n = 3 \bmod 4$,
	otherwise $B_n/x$ is symmetric. Thus, we may use
\[	\beta_k = \cases{\beta_{d-k},	&in the symmetric case;\cr
		         -\beta_{d-k},	&in the anti-symmetric case.\cr} \]
\item[]
Using these properties, the recurrences~(\ref{eq:AL3}--\ref{eq:AL4})
need only be applied for $k \le \max(1, \lfloor d/2 \rfloor)$.
\end{itemize}

\medskip
\leftline{\bf Example}
\medskip

\nopagebreak
Consider the case $n = 15$ as in Section~\ref{subsec:ex}.
We have $s = -1$, $d = \phi(15)/2 = 4$.  Thus
\[q_1 = q_2 = q_4 = \mu(15)\phi(1) = 1,\;
  q_3 = \mu(5)\phi(3) = -2,\]
\[r_1 = (1 \vert 15) = 1,\;
  r_2 = (2 \vert 15) = (2 \vert 3)(2 \vert 5) = 1,\;
  r_3 = (3 \vert 15) = 0,\;
  r_4 = (4 \vert 15) = 1 .\]
$q_3$, $q_4$, $r_3$, and $r_4$ are not required if we use symmetry.

\medskip						%
The initial conditions are $\alpha_0 = 2$ and $\beta_0 = 0$.
The recurrences~(\ref{eq:AL3}--\ref{eq:AL4}) give
\[\alpha_1 = (-15r_1\beta_0 - q_1\alpha_0)/2 = -1 ,\]
\[\beta_1  = (r_1\alpha_0 - q_1\beta_0)/2 = 1 ,\]
\[\alpha_2 =
  (-15r_2\beta_0 - 15r_1\beta_1 - q_2\alpha_0 - q_1\alpha_1)/4 = -4 , \]
\[\beta_2  = (r_2\alpha_0 + r_1\alpha_1 - q_2\beta_0 - q_1\beta_1)/4 = 0 .\]

Using symmetry of $A_n(x)$ and anti-symmetry of $B_n(x)/x$,
or continuing with the recurrences~(\ref{eq:AL3}--\ref{eq:AL4}),
we obtain
$\alpha_3 = \alpha_1 = -1$,
$\beta_3 = -\beta_1 = -1$,
$\alpha_4 = \alpha_0 = 2$,
$\beta_4 = -\beta_0 = 0$.
Thus $A_{15}(x) = 2x^4 - x^3 - 4x^2 - x + 2$ and $B_{15}(x) = x^3 - x$,
as expected.

\subsection{An algorithm for computing $C_n$ and $D_n$}
\label{subsec:AlgsCD}

Consider the computation of $C_n$ and $D_n$ for square-free $n > 1$.
Define $n'$, $s$, $s'$ and $L_n$ as in Section~\ref{subsec:Lucas},
and $d = \phi(n')/2$.  Thus ${\rm deg}\; L_n = 2d$, ${\rm deg}\; C_n = d$,
and ${\rm deg}\; D_n = d-1$.
From~(\ref{eq:A4.5}) it is enough to compute the coefficients $a_k$ of
$L_n(x)$.  In order to work over the integers, we define
\[	q_k = \cases{s'p_k/\sqrt{n}, 	&if $k$ is odd;\cr
		     p_k,		&if $k$ is even;\cr}	\]
where $p_k$ is the sum of $k$-th powers of roots of $L_n(x)$.
Thus, from~(\ref{eq:A5}),
\begin{equation}
  q_k = \cases{	(n \vert k),			 &if $k$ is odd;\cr
		\mu(n'/g_k')\phi(g_k')\cos\left((n-1)k\pi/4\right), &otherwise.\cr}
\label{eq:AL5}
\end{equation}
If
\[	C_n(x) = \sum_{j=0}^d \gamma_jx^{d-j} 	\]
and
\[	D_n(x) = \sum_{j=0}^{d-1} \delta_jx^{d-1-j} ,	\]
then, from~(\ref{eq:A4.5}),
\[	\gamma_k = a_{2k}	\]
and
\[	\delta_k = -s'a_{2k+1}/\sqrt{n} .	\]
In particular, $\gamma_0 = \delta_0 = 1$.
Using~(\ref{eq:NID}) we obtain
\begin{equation}
	\gamma_k = {1 \over 2k}\sum_{j=0}^{k-1}
		\left(nq_{2k-2j-1}\delta_j - q_{2k-2j}\gamma_j\right)
\label{eq:AL6}
\end{equation}
and
\begin{equation}
	\delta_k = {1 \over 2k+1}\left(\gamma_k + \sum_{j=0}^{k-1}
		\left(q_{2k+1-2j}\gamma_j - q_{2k-2j}\delta_j\right)\right)
\label{eq:AL7}
\end{equation}
for $k = 1, 2, \ldots$

We may use the fact that $C_n(x)$ and $D_n(x)$ are symmetric
to reduce the number of times the recurrences~(\ref{eq:AL6}--\ref{eq:AL7})
need to be applied.
An algorithm which incorporates this refinement is:
\medskip

\leftline{{\bf Algorithm L} (for Lucas)}

\nopagebreak
\begin{itemize}
\item[1.] Evaluate $q_k$ for $k = 1, \ldots, d$
	  using~(\ref{eq:AL5}).
\item[2.] Set $\gamma_0 \leftarrow 1$ and $\delta_0 \leftarrow 1$.
\item[3.] Evaluate $\gamma_k$ for $k = 1, \ldots, \lfloor d/2 \rfloor$
	  and $\delta_k$ for $k = 1, \ldots, \lfloor(d-1)/2\rfloor$
	  using~(\ref{eq:AL6}--\ref{eq:AL7}).
\item[4.] Evaluate $\gamma_k$ for $k = \lfloor d/2 \rfloor + 1, \ldots, d$
	  using $\gamma_k = \gamma_{d-k}$.
\item[5.] Evaluate $\delta_k$ for $k = \lfloor (d+1)/2 \rfloor, \ldots, d-1$
	  using $\delta_k = \delta_{d-1-k}$.
\end{itemize}

\medskip
\leftline{\bf Example}
\medskip

\nopagebreak
Consider the case $n = 15$ as in Section~\ref{subsec:ex}.
We have $n' = 2n = 30$, $s' = 1$, $d = \phi(30)/2 = 4$.
Thus
\[q_1 = (15 \vert 1) = 1 ,\]
\[q_2 = \mu(15)\phi(2)\cos(7\pi) = -1 ,\]
\[q_3 = (15 \vert 3) = 0 ,\]
\[q_4 = \mu(15)\phi(2)\cos(14\pi) = 1 .\]

The initial conditions are $\gamma_0 = \delta_0 = 1$.
The recurrences~(\ref{eq:AL6}--\ref{eq:AL7}) give
\[\gamma_1 = (15q_1\delta_0 - q_2\gamma_0)/2 = 8 ,\]
\[\delta_1 = (\gamma_1 + q_3\gamma_0 - q_2\delta_0)/3 = 3 ,\]
\[\gamma_2 =
  (15q_3\delta_0 + 15q_1\delta_1 - q_4\gamma_0 - q_2\gamma_1)/4 = 13 .\]

Using symmetry we obtain
$\gamma_3 = \gamma_1 = 8$,
$\gamma_4 = \gamma_0 = 1$,
$\delta_2 = \delta_1 = 3$,
and
$\delta_3 = \delta_0 = 1$.
Thus $C_{15}(x) = x^4 + 8x^3 + 13x^2 + 8x + 1$ and
$D_{15}(x) = x^3 + 3x^2 + 3x + 1$, as expected.

\subsection{Stevenhagen's algorithm}
\label{subsec:Stev}

Stevenhagen~\cite{Ste87} gives an algorithm for computing the polynomials
$C_n(x)$ and $D_n(x)$.  His algorithm depends on the application of the
Euclidean algorithm to two polynomials with integer coefficients 
and degree $O(n)$.
$C_n(x)$ and $D_n(x)$ may be computed as soon as a polynomial of
degree $\le \phi(n)/2$ is generated by the Euclidean algorithm.
Thus, the algorithm requires $O(n^2)$ arithmetic operations,
the same order\footnote{The complexity of
both algorithms can be reduced to $O(n (\log n)^2)$ arithmetic
operations by standard ``divide and conquer''
techniques~\cite{Aho74,Bre80}, but this is not of practical significance.}
as our Algorithm~L.

Unfortunately, Stevenhagen's algorithm suffers from a well-known
problem of the Euclidean algorithm~\cite{Knu81}~--
although the initial and final polynomials have small integer
coefficients, the intermediate results grow exponentially large.
When implemented in 32-bit integer arithmetic we found that
Stevenhagen's algorithm failed due to integer overflow
for $n = 35$.

Algorithm L does not suffer from this problem.
It is easy to see from the recurrences~(\ref{eq:AL6}--\ref{eq:AL7})
that intermediate results can grow only slightly larger than
the final coefficients $\gamma_k$ and $\delta_k$.
A straightforward implementation of
Algorithm~L can compute $C_n$ and $D_n$ for all square-free $n < 180$
without encountering integer overflow in 32-bit arithmetic.
When it does eventually occur,
overflow is easily detected because the division
by $2k$ in~(\ref{eq:AL6})
or by $2k+1$ in~(\ref{eq:AL7}) gives a non-integer result.

\section{Explicit expressions for $A_n, \ldots, D_n$}
\label{sec:GF}

We now use~(\ref{eq:GF1}) to give generating functions for
the coefficients of $A_n,\ldots,D_n$.
The generating functions can be used to evaluate the
coefficients of $A_n(x),\ldots,D_n(x)$ in $O(n \log n)$ arithmetic
operations, via the fast power series algorithms of Section 5 of
Brent and Kung~\cite{Bre78}.		%
Also, where the generating functions
converge, they give explicit formulas which can be used to
compute $A_n(x),\ldots,D_n(x)$ at particular arguments $x$.
However, it is often more efficient to compute
the coefficients of the polynomials
by the algorithms of Section~\ref{sec:Algs}
and then evaluate the polynomials by Horner's rule.

The generating functions may be written in terms of certain
analytic functions $f_n$ and $g_n$, which we now define.

\subsection{The analytic functions $f_n$ and $g_n$}
\label{subsec:fg}
For odd square-free $n > 1$ and $\vert x \vert \le 1$, define
\begin{equation}
f_n(x) = \sum_{j=1}^\infty{(j \vert n){x^j \over j}} . \label{eq:f}
\end{equation}
Similarly, for square-free $n > 1$ and $\vert x \vert \le 1$, define
\begin{equation}
g_n(x) = \sum_{j=0}^\infty{(n \;\vert\; {2j+1})
			{x^{2j+1} \over {2j+1}}} . 		\label{eq:g}
\end{equation}
Observe that $g_n(x)$ is an odd function, so $g_n(-x) = -g_n(x)$.

It follows from~(\ref{eq:GoverG}) and~(\ref{eq:GL1}) below that
$\exp(\sqrt{sn}f_n(x))$ and $\exp(2\sqrt{n}g_n(x))$
are rational functions with zeros and poles at certain roots of unity.
From these representations it follows that the analytic continuations
outside the unit circle are given by
\begin{equation}
\cases {f_n(x) = f_n(1/x), 		&if $n = 1 \bmod 4$;\cr
        f_n(x) + f_n(1/x) = 2f_n(1),    &if $n = 3 \bmod 4$;\cr}
							\label{eq:fcont}
\end{equation}
and
\begin{equation}
	g_n(x) = g_n(1/x) .				\label{eq:gcont}
\end{equation}

The functions $f_n(x)$ and $g_n(x)$ are closely related.
For example, taking the odd terms in the sum~(\ref{eq:f}) and using
the law of quadratic reciprocity, we obtain
\[ f_n(x) - f_n(-x) = 2g_n(x\sqrt{s})/\sqrt{s} .	\]
Such identities are a consequence of
relationships between the polynomials~$G_n(x)$ and~$L_n(x)$.

$f_n(1)$ is related to the
class number $h(D)$ of the quadratic field $Q[\sqrt{D}]$
with discriminant $D = sn$.
In the notation of Davenport~\cite{Dav80},	%
$f_n(1) = L_{-1}(1) = L(1) = L(1,\chi)$,
where $\chi(j) = (j \vert n)$ is the real,
nonprincipal Dirichlet character appearing in~(\ref{eq:f}).
Known results~\cite{Dav80,Lan27,Lan90,Was82}
in the case $n = 3 \bmod 4$ (so $D = -n$) are
\begin{equation}
 f_n(1) = %
	{-\pi \over {n^{3/2}}}\sum_{j=1}^{n-1} (j \vert n)j =
	{\pi \over (2 - (2 \vert n))\sqrt{n}}
		\sum_{j=1}^{(n-1)/2} (j \vert n) =
	{2\pi \over w\sqrt{n}}h(-n) > 0 . 		\label{eq:hn3}
\end{equation}
Here
\[w = \cases{6, &if $n = 3$;\cr
	     2, &if $n = 3 \bmod 4$, $n > 3$\cr}	\]
is the number of roots of unity in $Q[\sqrt{D}]$.
Since $h(-n)$ is an integer and $h(-3) = 1$, we have
\begin{equation}
\exp(2i\sqrt{n}f_n(1)) = \cases{(-1 + \sqrt{-3})/2, &if $n = 3$;\cr
		1, &if $n = 3 \bmod 4$, $n > 3.$\cr}		\label{eq:exp1}
\end{equation}

In the case $n = 1 \bmod 4$, we have $D = n$, and
\begin{equation}
f_n(1) = %
	{{\ln \varepsilon} \over \sqrt{n}}h(n) .	\label{eq:hn1}
\end{equation}
Here $\varepsilon$ is the ``fundamental unit'',
i.e.~$\varepsilon = (\vert u \vert + \sqrt{n}\vert v \vert)/2$,
where $(u,v)$ is a minimal nontrivial solution of
$u^2 - nv^2 = 4$.
For example, if $n = 5$ then $\varepsilon = (3 + \sqrt{5})/2$, $h(5) = 1$ and
$f_5(1) = (\ln \varepsilon)/\sqrt{5} = 0.4304\ldots$

\subsection{The polynomials $A_n$ and $B_n$}
\label{subsec:AB}

Let $n > 3$ be odd and square-free. We exclude $n = 3$ to avoid the special
case in~(\ref{eq:exp1}), but the results apply with minor modifications
when $n = 3$.  Let $s$, $G_n$, ${\tilde G}_n$ be as in
Section~\ref{subsec:Gauss}, and $d = \phi(n)/2$.
Recall that
\begin{equation}
  G_n(x) = \prod_{\scriptstyle 0 < j < n \atop
	\scriptstyle  (j \vert n)=1} (x - \zeta^j) 		\label{eq:G2b}
\end{equation}
and
\begin{equation}
  {\tilde G}_n(x) = \prod_{\scriptstyle 0 < j < n \atop
	\scriptstyle  (j \vert n)=-1} (x - \zeta^j) .		\label{eq:G2c}
\end{equation}
From~(\ref{eq:GF1}) and~(\ref{eq:G6}), we have
\begin{equation}
{{\tilde G}_n(1/x) / G_n(1/x)} = \exp(\sqrt{sn}f_n(x)) .
							     \label{eq:GoverG}
\end{equation}
Also, from~(\ref{eq:G2b}),
\begin{equation}
(-x)^dG_n(1/x) = \prod_{\scriptstyle 0 < j < n \atop
	\scriptstyle  (j \vert n)=1} (\zeta^{j}x - 1) =
		 \prod_{\scriptstyle 0 < j < n \atop
	\scriptstyle  (j \vert n)=1} \zeta^j (x - \zeta^{-j}) ,	\label{eq:G2d}
\end{equation}
so
\begin{equation}
{G_n(1/x) / {\tilde G}_n(1/x)} = \zeta^{\sigma} \prod_{j=1}^{n-1}
	(x - \zeta^{-j})^{(j \vert n)} ,			\label{eq:G2e}
\end{equation}
where
\begin{equation}
\sigma = \sum_{j=1}^{n-1}(j \vert n)j .		\label{eq:sigma}
\end{equation}
If $n = 1 \bmod 4$ then by grouping the terms for $j$ and $n-j$ ($j < n/2$)
in~(\ref{eq:sigma}) we have $n \vert \sigma$. If $n = 3 \bmod 4$ then,
from~(\ref{eq:hn3}), we have $\sigma = -nh(-n)$ so again $n \vert \sigma$.
Thus, in both cases $\zeta^{\sigma} = 1$,
and from~(\ref{eq:G2e}) we have
\begin{equation}
{G_n(1/x) / {\tilde G}_n(1/x)} =
	\cases{ {G_n(x) / {\tilde G}_n(x)}, &if $n = 1 \bmod 4$;\cr
	        {{\tilde G}_n(x) / G_n(x)}, &if $n = 3 \bmod 4$.\cr}
\end{equation}
It follows from~(\ref{eq:GoverG}) that
\begin{equation}
{{\tilde G}_n(x) / G_n(x)} = \exp(s\sqrt{sn}f_n(x)) .
							     \label{eq:GoverG2}
\end{equation}
We see from~(\ref{eq:GoverG}) or~(\ref{eq:GoverG2}) that, as claimed above,
$\exp(\sqrt{sn}f_n(x))$ is a rational function. It has zeros at
$\zeta^j$, $(j \vert n) = -s$;
and poles at
$\zeta^j$, $(j \vert n) = +s$.
From~(\ref{eq:G2.6}) and~(\ref{eq:GoverG2}), taking a square root,
we obtain
\begin{equation}	%
G_n(x) = \sqrt{\Phi_n(x)}\exp\left({{-s\sqrt{sn}}\over{2}}f_n(x)\right)
								\label{eq:GenG}
\end{equation}
If~(\ref{eq:GenG}) is interpreted as a generating function for $G_n(x)$
then $\sqrt{\Phi_n(x)}$ and $f_n(x)$
should be interpreted as power series in $x$,
and the correct sign of the square root is positive.  On the other hand,
if~(\ref{eq:GenG}) is regarded as an exact expression for $G_n(x)$,
then the sign of the square root
is positive for real $x$, because $G_n(x)$ and
$\Phi_n(x)$ have no real roots, and the exponential never vanishes,
so a change in sign would contradict the continuity of $G_n(x)$.
An extension of this argument shows that the same branch of the square root
must be taken in any simply-connected, closed region which does
not contain any of the zeros of $\Phi_n(x)$.
(We omit similar comments below.)

From~(\ref{eq:GenG}) we easily deduce the corresponding
expressions for $A_n(x) = G_n(x) + {\tilde G}_n(x)$
and $B_n(x) = ({\tilde G}_n(x) - G_n(x))/\sqrt{sn}$. We state the results
as a Theorem:

\newtheorem{thmab}{Theorem}
\begin{thmab}
For odd, square-free $n > 3$, the polynomials $A_n(x)$ and $B_n(x)$
occurring in Gauss's identity~{\rm{(\ref{eq:G1})}} are
\begin{equation}
A_n(x) = 2\sqrt{\Phi_n(x)}\cosh\left({{\sqrt{sn}}\over{2}}f_n(x)\right)
						\label{eq:GF2}
\end{equation}
and
\begin{equation}
B_n(x) = 2\sqrt{\Phi_n(x) \over sn}\sinh\left({{\sqrt{sn}}\over{2}}
		f_n(x)\right) .
						\label{eq:GF3}
\end{equation}
\label{thm:Tab}
\end{thmab}

\leftline{\bf Remark}
\medskip
\nopagebreak
If $n = 3 \bmod 4$, so $s = -1$, then
it is natural to replace $\cosh(iz)$ by $\cos(z)$ in~(\ref{eq:GF2})
and $\sinh(iz)$ by $i\sin(z)$ in~(\ref{eq:GF3}), giving
\begin{equation}
A_n(x) = 2\sqrt{\Phi_n(x)}\cos\left({{\sqrt{n}}\over{2}}f_n(x)\right)
						\label{eq:GF2b}
\end{equation}
and
\begin{equation}
B_n(x) = 2\sqrt{\Phi_n(x) \over n}\sin\left({{\sqrt{n}}\over{2}}
		f_n(x)\right) .
						\label{eq:GF3b}
\end{equation}

\medskip
\leftline{\bf Example}
\medskip

\nopagebreak
Consider the case $n = 15$.  We expand the right side of~(\ref{eq:GF2b})
as a power series in $x$, keeping enough terms to find $A_{15}(x)$
without using symmetry.
From Section~\ref{subsec:ex} we have
\[\Phi_{15}(x) = 1 - x + x^3 - x^4 + x^5 - x^7 + x^8, \]
so
\[\sqrt{\Phi_{15}(x)} = 1 - {x \over 2} - {x^2 \over 8} + {7x^3 \over 16}
	- {37x^4 \over 128} + \cdots\]
Also,
\[f_{15}(x) = x + {x^2 \over 2} + {x^4 \over 4} - {x^7 \over 7} + \cdots \]
so
\[\cos\left({\sqrt{15} \over 2}f_{15}(x)\right) =
  1 - {15x^2 \over 8} - {15x^3 \over 8} + {15x^4\over 128} + \cdots \]
and
\[2\sqrt{\Phi_{15}(x)}\cos\left({\sqrt{15}\over 2}f_{15}(x)\right) =
  2 - x - 4x^2 - x^3 + 2x^4 + \cdots,\]
which is $A_{15}(x) + O(x^5)$ as expected.   
We can ignore the $O(x^5)$ term (which in fact vanishes) since we know that
${\rm deg}\; A_{15}(x) = \phi(15)/2 = 4$.

A reader who attempts similar computations for larger $n$ will soon be
convinced that Algorithm~D of Section~\ref{subsec:AlgsAB} is more
convenient, if only because all intermediate results are integers
and there is an easy check on the accuracy of the inner product
accumulations.

Using~(\ref{eq:hn3}--\ref{eq:hn1}) and the fact
(an immediate consequence of~(\ref{eq:phirec})) that
\[ \Phi_n(1) = \cases{n, &if $n$ is prime;\cr
		     1, &for other square-free $n > 1$;\cr}		\]
we can verify that~(\ref{eq:GoverG}) and~(\ref{eq:GoverG2}--\ref{eq:GF3b})
give sensible results in the case $x = 1$.
For example, we have $h(-15) = 2$, so $f_{15}(1) = 2\pi/\sqrt{15}$
and $\cos({\sqrt{15} \over 2}f_{15}(1)) = \cos(\pi) = -1$.
Thus~(\ref{eq:GF2b}) gives the correct value $-2$ of $A_{15}(1)$.

\subsection{The polynomials $C_n$ and $D_n$}
\label{subsec:CD}

We now consider the analogues of~(\ref{eq:GF2}--\ref{eq:GF3}) for
the Lucas polynomials $C_n$ and $D_n$. The argument is similar to that
of Section~\ref{subsec:AB}, but simpler because the polynomial
$L_n(x)$ is symmetric, which leads to the
simple functional equation~(\ref{eq:gcont}) for $g_n(x)$.

Assume that $n > 1$ is square-free, and adopt the notation of
Section~\ref{subsec:Lucas}.  Using~(\ref{eq:A5}) in the
generating function~(\ref{eq:GF1}), we have
\begin{equation}
{\tilde L}_n(1/x) / L_n(1/x) = \exp\left(2s'\sqrt{n}g_n(x)\right) .
								 \label{eq:GL1}
\end{equation}
This shows that $\exp(2\sqrt{n}g_n(x))$ is a rational function.
From
\[{\tilde L}_n(1/x) / L_n(1/x) = {\tilde L}_n(x) / L_n(x)	\]
we deduce the functional equation~(\ref{eq:gcont}), which gives the analytic
continuation of $g_n(x)$ outside the unit circle. We may also
write~(\ref{eq:GL1}) more simply as
\begin{equation}
{\tilde L}_n(x) / L_n(x) = \exp\left(2s'\sqrt{n}g_n(x)\right) .
								 \label{eq:GL2}
\end{equation}
From~(\ref{eq:A4.7}) and~(\ref{eq:GL2}), taking a square root, we obtain
\begin{equation}
L_n(x) = \sqrt{F_n(x^2)}\exp\left(-s'\sqrt{n}g_n(x)\right) . 	\label{eq:GL3}
\end{equation}

\newtheorem{thmcd}[thmab]{Theorem}
\begin{thmcd}
Let $n > 1$ be square-free.
The Aurifeuillian factors
$F_n^{\pm}(x) = C_n(x) \pm\sqrt{nx}D_n(x)$ of $F_n(x)$ are given by
\begin{equation}
F_n^{\pm}(x) = \sqrt{F_n(x)}\exp\left(\pm\sqrt{n}g_n(\sqrt{x})\right) .
							\label{eq:Aufacs}
\end{equation}
Also,
\begin{equation}
C_n(x) = \sqrt{F_n(x)}\cosh\left(\sqrt{n}g_n(\sqrt{x})\right)	\label{eq:GL5}
\end{equation}
and
\begin{equation}
D_n(x) = \sqrt{F_n(x) \over nx}\sinh\left(\sqrt{n}g_n(\sqrt{x})\right) .
								\label{eq:GL6}
\end{equation}
\label{thm:Tcd}
\end{thmcd}

\leftline{\bf Proof}
\medskip
\nopagebreak
Recall that the Aurifeuillian factors
$F_n^{\pm}(x) = C_n(x) \pm\sqrt{nx}D_n(x)$ of $F_n(x)$ are
$L_n(\pm\sqrt{x})$.
Thus~(\ref{eq:Aufacs}) follows from~(\ref{eq:GL3}).
Since
\[	L_n(x) + L_n(-x) = 2C_n(x^2)	\]
and
\[	L_n(x) - L_n(-x) = -2s'x\sqrt{n}D_n(x^2),	\]
we easily deduce~(\ref{eq:GL5}--\ref{eq:GL6}).				\QED

\subsection{Application to integer factorization}
\label{subsec:intfac}

In this section we illustrate how the results
of Sections~\ref{subsec:AlgsCD} and~\ref{subsec:CD} can be used to obtain
factors of integers of the form $a^n \pm b^n$. Our examples are
for illustrative purposes, so are small enough to be verified by hand.
Many larger examples can be found in~\cite{Bre92,Bri83}.

As usual, $n > 1$ is a square-free integer.
Recall the definition of $F_n(x)$ from Section~\ref{subsec:no}.
Note that the polynomial $F_n(x)$ is a factor of $x^n \pm 1$
(where the sign is ``$-$'' if $n = 1 \bmod 4$, and ``$+$'' otherwise).

If $x$ has the form $m^2 n$, where $m$ is a positive integer,
then $\sqrt{nx} = mn$ is an integer, and the Aurifeuillian factors
$F_n^{\pm}(x) = C_n(x) \pm mnD_n(x)$
give integer factors of $F_n(x)$,
and hence of $x^n \pm 1 = m^{2n}n^n \pm 1$.
For example, if $m = n^k$, we obtain factors of
$n^{(2k+1)n} \pm 1$.

More generally, if $m = p/q$ is rational, we obtain rational factors
of $x^n \pm 1 = p^{2n}q^{-2n}n^n \pm 1$,
and thus integer factors of
$p^{2n}n^n \pm q^{2n}$.
We consider one example later, but for the moment
we continue to assume that $m$ is an integer.

Before giving numerical examples, we state explicitly how the results of
Section~\ref{subsec:CD} can be used to compute $F_n^{\pm}(m^2 n)$.

\newtheorem{thma}[thmab]{Theorem}
\begin{thma}
Let $m$, $n$ be positive integers, $n > 1$ square-free,
$x = m^2n$, and
$\lambda = \phi(2n)/2$ .
Then the Aurifeuillian factors $F_n^{\pm}(x)$ of $F_n(x)$ are given by
\[ F_n^{-}(x) = \lfloor{\hat F} + \;1/2 \rfloor \]
and
\[ F_n^{+}(x) = F_n(x)/F_n^{-}(x) ,	\]
where
\[ {\hat F} = \sqrt{F_n(x)}\exp\left(
	-\;{1 \over m}\sum_{j=0}^{\lambda-1}
	{(n\;\vert\;2j+1) \over (2j+1)x^j}\right) . 	\]
\label{thm:Ta}
\end{thma}

\leftline{\bf Proof}
\medskip
\nopagebreak
From~(\ref{eq:Aufacs}),
using the functional equation~(\ref{eq:gcont})
and the power series~(\ref{eq:g}) for $g_n(1/\sqrt{x})$,
we have
\[ F_n^{-}(x) = \sqrt{F_n(x)}\exp\left(
	-\;{1 \over m}\sum_{j=0}^{\infty}
	{(n\;\vert\;2j+1) \over (2j+1)x^j}\right) . 	\]
Thus, we only have to show that
the error incurred by truncating
the power series after $\lambda$ terms is less than $1/2$ in absolute
value, i.e.~that
$\vert {\hat F} - F_n^{-} \vert < 1/2$.

If $x < 5$ we must have $m = 1$, $n = 2$ or $3$, $x = n$, $\lambda = 1$.
In both cases $F_n^{-}(n) = 1$ and it is easy to verify that
$1/2 < {\hat F} < 3/2$.  Thus, from now on we assume that $x \ge 5$.

Let
\[ t ={1 \over m}\sum_{j=\lambda}^{\infty}
	{(n\;\vert\;2j+1) \over (2j+1)x^j} .	\]
Since $m \ge 1$, $\lambda \ge 1$, $x \ge 5$,
and the Jacobi symbol is 0 or $\pm 1$, we have
\[ \vert t \vert \le \sum_{j=\lambda}^{\infty}{x^{-j} \over (2j+1)}
	\le \left({1 \over 3} + {1 \over 5 \cdot 5} +
	{1 \over 7 \cdot 5^2} + \cdots\right)x^{-\lambda}	\]
so
\begin{equation}
        \vert t \vert
	\le 5\left(\sqrt{5}\;\ln\left({1 + \sqrt{5} \over 2}\right) - 1\right)
		x^{-\lambda} < 0.3802x^{-\lambda} .
							\label{eq:crude}
\end{equation}
In particular, $\vert t \vert < 0.3802/x < 0.08$, so
\begin{equation}
  \vert \exp(t) - 1 \vert < {\vert t \vert \over 1 - \vert t \vert/2}
	< {0.3802x^{-\lambda} \over 0.96}
	< 0.4x^{-\lambda} .					\label{eq:ina}
\end{equation}
Now ${\hat F} = F_n^{-}\exp(t)$, so
\begin{equation}
 \vert {\hat F} - F_n^{-} \vert = F_n^{-}\vert \exp(t) - 1 \vert .
								\label{eq:inb}
\end{equation}
Applying Lemma~\ref{lemma:Lb} of Section~\ref{subsec:cyclo} with $R = 5$ gives
\[ F_n(x) < \exp(1/4)x^{\phi(2n)} ,	\]
but $\phi(2n) = 2\lambda$ and $F_n^{-}$ is the smaller factor of $F_n$, so
\begin{equation}
 F_n^{-}(x) \le \sqrt{F_n(x)} < \exp(1/8)x^\lambda .		\label{eq:inc}
\end{equation}
From~(\ref{eq:ina}--\ref{eq:inc}) we finally obtain the bound
\[ \vert {\hat F} - F_n^{-} \vert < 0.4\exp(1/8) < 0.5 ,	\]
which completes the proof.					\QED
\medskip

Since the argument of the exponential in Theorem~\ref{thm:Ta} 
is $-1/m + O(1/n)$
as $n \rightarrow \infty$,
we have the following result, which sheds some light on the
ratio of the Aurifeuillian factors. (The Corollary strictly follows
from Theorem~\ref{thm:Ta} only if $m$ is an integer, but from the proof
of Theorem~\ref{thm:Ta} it is clear that the Corollary is also valid
for rational $m$.)

\newtheorem{cora}{Corollary}
\begin{cora}
Let $x = m^2n$ where $m > 0$ is rational
$n > 1$ is an integer.
Consider $m$ fixed as $n \rightarrow \infty$ through square-free values.
Then the Aurifeuillian factors $F_n^{\pm}(x)$ of $F_n(x)$ satisfy
\[ F_n^{+}(x)/F_n^{-}(x) = \exp(2/m) + O(1/n) .	\]
\label{cor:Ca}
\end{cora}

\leftline{\bf Examples}
\medskip

\nopagebreak
{\bf 1.}\ To start with a simple example, consider
$n = 2$, $m = 2$, so $x = m^2n = 8$
and $\lambda = 1$.
In Theorem~\ref{thm:Ta} we have $F_2(x) = x^2 + 1 = 65$,
\[ {\hat F} = \sqrt{F_2(x)}\exp(-1/m) =
	\sqrt{65}\exp(-1/2) = 4.89\cdots	\]
and rounding to the nearest integer gives the factor~5 of 65.
\medskip

{\bf 2.}\ A similar but less trivial example is
$n = 2$, $m = 2^5$, $x = m^2n = 2^{11}$.
Here $F_2(x) = x^2 + 1 = 2^{22}+1$ and
\[ {\hat F} = \sqrt{F_2(x)}\exp(-1/m) =
	\sqrt{2^{22}+1}\exp(-2^{-5}) = 1984.98\cdots	\]
which on rounding gives the factor $F_2^{-} = 1985$ of $2^{22}+1$.
By division we find $F_2^{+} = 2113$, so the complete factorization
is $2^{22}+1 = 5 \cdot 397 \cdot 2113$ .
\medskip

{\bf 3.}\ Now consider $n = 5$, $m = 3$, so $x = m^2n = 45$ and
$\lambda = \phi(10)/2 = 2$.
In this case
\[ 	F_5(x) = \Phi_5(x) = (x^5 - 1)/(x - 1) = 4193821 ,	\]
\[ {\hat F} = \sqrt{F_5(x)}\exp\left(-{1 \over m} + {1 \over 3m^3n}\right) =
	\sqrt{4193821}\exp(-134/405) = 1470.99924\cdots	\]
and rounding to the nearest integer gives the factor~1471 of $F_5(x)$.
By division we obtain the other factor~2851. In this example, but not
in general, the Aurifeuillian factors are prime.
\medskip

{\bf 4.}\ Now consider a composite $n$, say $n = 15$.  To keep the arithmetic easy
we take $m = 1$, so $x = 15$ and
\[ F_{15}(x) = \Phi_{15}(-x)
	= x^8 + x^7 - x^5 - x^4 - x^3 + x + 1
	= 2732936641 .	\]
The example in Section~\ref{subsec:AlgsCD} shows how we can
use Algorithm~L to compute
\[ C_{15}(x) = x^4 + 8x^3 + 13x^2 + 8x + 1 	\]
and
\[ D_{15}(x) = x^3 + 3x^2 + 3x + 1 .		\]
Evaluating the polynomials, we find
\[ C_{15}(15) = 80671, \;\; D_{15}(15) = 4096 ,	\]
so the Aurifeuillian factors of $F_{15}(15)$ are
$80671 \pm 15 \cdot 4096$. This gives 19231 and 142111, which
can easily be verified to be factors of $15^{15}+1$.
The ``algebraic'' factors $(15^3 + 1)/16$ and $(15^5 + 1)/16$ allow us to
complete the factorization:
\[ 15^{15}+1 =
   2^4 \cdot 31 \cdot 211 \cdot 1531 \cdot 19231 \cdot 142111	\]

Alternatively, instead of evaluating $C_{15}(15)$ and $D_{15}(15)$,
we can use Theorem~\ref{thm:Ta}.
We have $\lambda = \phi(15)/2 = 4$.
Since $(15 \vert 3) = (15 \vert 5) = 0$ and $(15 \vert 7) = 1$,
the computation gives
\[ {\hat F} = \sqrt{F_{15}(x)}\exp\left(-1 - {1 \over 7n^3}\right) =
	19231.00217\cdots	\]
and rounding to the nearest integer gives the
factor~19231 of $F_{15}(15)$.
\medskip

{\bf 5.}\ To conclude, we give an example where $m = p/q$ is rational but not
an integer.  Consider $n = 7$, $p = 2$, $q = 5$,
so $x = p^2 n/q^2 = 28/25$
and $\sqrt{nx} = pn/q = 14/5$.
We have
\[ F_7(x) = \Phi_7(-x) = {x^7 + 1 \over x + 1}
	= {28^7 + 25^7 \over 53 \cdot 25^6} 	\]
Theorem~\ref{thm:Ta} is not applicable.
Because $x$ is close to 1, the
series for $g_7(1/\sqrt{x})$ converges rather slowly,
and we need to take at least 35 terms			%
to obtain sufficient accuracy.
However, using Algorithm~L we easily find that
\[ C_7(x) = x^3 + 3x^2 + 3x + 1, \;\; D_7(x) = x^2 + x + 1,	\]
so Horner's rule gives
\[ C_7(x) = {148877/5^6}, \;\; D_7(x) = {2109/5^4} , \]
and $C_7(x) \pm \sqrt{nx}D_7(x)$ gives
$296507/5^6$ and $1247/5^6$.
Thus, we have obtained factors 296507 and 1247 of $(25^7 + 28^7)/53$.
The larger factor is prime, so it is easy to deduce the complete
factorization 
\[25^7 + 28^7 = 29 \cdot 43 \cdot 53 \cdot 296507	\]

\medskip
\leftline{\bf Acknowledgements}
\medskip

\nopagebreak
Thanks are due to Brendan McKay for suggesting the use of
the generating function~(\ref{eq:GF1}),
to Emma Lehmer and an anonymous referee
for helpful comments on the exposition,
to Hans Riesel for his kind assistance with the solution of
exercise A6.2 of~\cite{Rie85},
and to Sam Wagstaff and Hugh Williams for providing copies of
several references which were difficult to find in Australia.


\begin{thebibliography}{99}

\bibitem{Aho74}
A.\ V.\ Aho, J.\ E.\ Hopcroft and J.\ D.\ Ullman,
{\em The Design and Analysis of Computer Algorithms},
Addison-Wesley, Menlo Park, Calif., 1974, Chapter 8.

\bibitem{Aur77}
A.\ Aurifeuille and H.\ Le Lasseur, see~\cite{Luc77}, page 276,
or~\cite{Luc78a}, page 785.

\bibitem{Bee51}
N.\ G.\ W.\ H.\ Beeger,
``On a new quadratic form for certain cyclotomic polynomials'',
{\em Nieuw Arch.\ Wisk.} (2), 23 (1951), 249-252.

\bibitem{Bre80}
R.\ P.\ Brent, F.\ G.\ Gustavson and D.\ Y.\ Y.\ Yun,
``Fast solution of Toeplitz systems of equations and computation of
Pad\'e approximants'',
{\em J.\ Algorithms} 1 (1980), 259-295.

\bibitem{Bre78}
R.\ P.\ Brent and H.\ T.\ Kung,
``Fast algorithms for manipulating formal power series'',
{\em J.~ACM} 25 (1978), 581-595.

\bibitem{Bre92}
R.\ P.\ Brent and H.\ J.\ J.\ te Riele,
{\em Factorizations of $a^n \pm 1$, $13 \le a < 100$},
Report NM-R9212, Department of Numerical Mathematics,
Centrum voor Wiskunde en Informatica,
Amsterdam, June 1992. %

\bibitem{Bri83}
J.\ Brillhart, D.\ H.\ Lehmer, J.\ L.\ Selfridge, B.\ Tuckerman and
S.\ S.\ Wagstaff, Jr.,
{\em Factorizations of $b^n \pm 1$,
$b = 2, 3, 5, 6, 7, 10, 11, 12$ up to high powers}, second edition,
American Mathematical Society, Providence, Rhode Island, 1988.


\bibitem{Cun15}
A.\ J.\ C.\ Cunningham,
``Factorisation of $N = y^y \mp 1$ and $x^{xy} \mp y^{xy}$'',
{\em Messenger of Math.} (2), 45 (1915), 49-75.

\bibitem{Cun25}		%
A.\ J.\ C.\ Cunningham and H.\ J.\ Woodall,
{\em Factorisation of $y^n \mp 1$, $y = 2, 3, 5, 6, 7, 10, 11, 12$
up to high powers $(n)$}, Hodgson, London, 1925.

\bibitem{Dav80}
H.\ Davenport,
{\em Multiplicative Number Theory},
second edition (revised by H.\ L.\ Montgomery),
Springer-Verlag, New York, 1980.

\bibitem{Dir94}
P.\ G.\ Lejeune Dirichlet,
{\em Vorlesungen \"uber Zahlentheorie}, fourth edition,
Friedr.\ Vieweg \& Sohn, Braunschweig, 1894, Chapter~5 and Supplement~7.

\bibitem{Gau01}
C.\ F.\ Gauss,
{\em Disquisitiones Arithmetic{\ae}}, G.\ Fleischer, Leipzig, 1801,
Art.~356-357.
Reprinted in {\em Carl Friedrich Gauss Werke}, Band 1,
Georg Olms Verlag, Hildesheim, 1981.

\bibitem{Har60}
G.\ H.\ Hardy and E.\ M.\ Wright,
{\em An Introduction to the Theory of Numbers},
fifth edition, Clarendon Press, Oxford, 1984, Ch.~16.

\bibitem{Knu81}
D.~E.~Knuth,
{\em The Art of Computer Programming,
Volume~2: Seminumerical Algorithms} (second edition),
Addison-Wesley, Menlo Park, 1981, Chapter 3.

\bibitem{Kra22}
M.\ Kraitchik,
``D\'ecomposition de $a^n \pm b^n$ en facteurs dans le cas o\`u
$nab$ est un carr\'e parfait avec une table des d\'ecompositions
num\'eriques pour toutes les valeurs de $a$ et $b$ inf\'erieures \`a 100'',
Gauthiers-Villars,	%
Paris, 1922.

\bibitem{Kra24}
M.\ Kraitchik,
{\em Recherches sur la Th\'eorie des Nombres}, Volume~1,
Gauthiers-Villars,	%
Paris, 1924.

\bibitem{Kra52}
M.\ Kraitchik,
{\em Introduction \`a la Th\'eorie des Nombres},
Gauthiers-Villars, Paris, 1952.

\bibitem{Lan27}
Edmund Landau,
{\em Vorlesungen \"uber Zahlentheorie}, Band 1(1):
{\em Aus der elementaren und additiven Zahlentheorie},
Leipzig, 1927. English translation: {\em Elementary Number Theory},
Chelsea, New York, 1958.


\bibitem{Lan90}
Serge Lang,
{\em Cyclotomic Fields I and II},
combined second edition, Graduate Texts in Mathematics 126,
Springer-Verlag, New York, 1990.

\bibitem{Luc77}
E.\ Lucas, ``Th\'eor\`emes d'arithm\'etique'',
{\em Atti.\ R.\ Acad.\ Sc.\ Torino} 13 (1877-8), 271-284.

\bibitem{Luc78a}
E.\ Lucas, ``Sur la s\'erie r\'ecurrente de Fermat'',
{\em Bull.\ Bibl.\ Storia Sc.\ Mat.\ e Fis.} 11 (1878), 783-789.

\bibitem{Luc78b}
E.\ Lucas, ``Sur les formules de Cauchy et de Lejeune-Dirichlet'',
{\em Ass.\ Fran\c caise
pour l'Avanc.\ des Sci., Comptes Rendus} 7 (1878), 164-173.

\bibitem{Ram18}		%
S.\ Ramanujan, ``On certain trigonometrical sums and their applications
in the theory of numbers'',
{\em Trans.\ Cambridge Philos.\ Soc.\ }22, 13 (1918), 259-276.
Reprinted in {\em Collected Papers of Srinivasa Ramanujan}
(edited by G.\ H.\ Hardy, P.\ V.\ Seshu Aiyar and B.\ M.\ Wilson),
Cambridge Univ.\ Press, 1927.	%

\bibitem{Rie85}
Hans Riesel,
{\em Prime Numbers and Computer Methods for Factorization},
Birkh\"auser, Boston, 1985.

\bibitem{Rio78}
John Riordan,
{\em An Introduction to Combinatorial Analysis},
Princeton University Press, New Jersey, 1978, Ch.~2, exercise~27.

\bibitem{Sch62}
A.\ Schinzel,
``On the primitive prime factors of $a^n - b^n$'',
{\em Proc.\ Cambridge Philos.\ Soc.\ }58 (1962), 555-562.
MR~26\#1280

\bibitem{Ste87}
Peter Stevenhagen,
``On Aurifeuillian factorizations'',
{\em Nederl.\ Akad.\ Wetensch.\ Indag.\ Math.} 49 (1987), 451-468.
MR~89a:11015	%

\bibitem{Tur52}
H.\ W.\ Turnbull,
{\em Theory of Equations} (fifth edition),
Oliver and Boyd, Edinburgh, 1952, Sec.~32.

\bibitem{Wae53}
B.\ L.\ van der Waerden,
{\em Algebra}, Vol.\ 1 (English translation by Fred Blum, fifth edition),
Frederick Ungar, New York, 1953, Ch.~7.

\bibitem{Was82}
L.\ C.\ Washington,
{\em Introduction to Cyclotomic Fields},
Graduate Texts in Mathematics 83, Springer-Verlag, New York, 1982.

\end{thebibliography}
\end{document}